\documentclass[12pt]{amsart}
\usepackage{amssymb}
\usepackage{amsfonts}
\usepackage{amsmath}
\usepackage{hyperref}
\hypersetup{hidelinks}

\newcommand{\func}[1]{\operatorname{#1}}

\theoremstyle{plain}
\newtheorem{acknowledgement}{Acknowledgement}

\newtheorem{definition}{Definition}

\newtheorem{remark}{Remark}

\numberwithin{equation}{section}

\title{The isoperimetric inequality for the Ky Fan norm}
\author{Lu\'{\i}s Daniel Abreu}
\address{NuHAG, Faculty of Mathematics, University of Vienna,
Oskar-Morgenstern-Platz 1, A-1090, Vienna, Austria}
\email{abreuluisdaniel@gmail.com}
\date{}

\begin{document}
\begin{abstract}
We show that, among all measurable sets $\Omega \subset \mathbb{C}$ with
finite area $s$, the disc of area $s$ maximizes the Ky Fan norm, which is
defined as the sum of the first $N$ eigenvalues of the Toeplitz operator
with symbol $\mathbf{1}_{\Omega }$ on the Fock space. For $N=1$\ this
reduces to Nicola-Tilli's celebrated Faber--Krahn inequality and for general 
$N$ the result was conjectured by Nicola, Riccardi and Tilli, who proved it
for radial sets. The proof combines Ky Fan's maximum principle with methods
from quantum information theory based on Fock rearrangements of density
operators. As a by-product, we obtain isoperimetric inequalities for the
Schatten sums in the range $0<p\leq \infty $.
\end{abstract}
\maketitle

\section{Introduction}

Consider the Bargmann--Fock space of entire functions, defined as%
\begin{equation*}
\mathcal{F}^{2}(\mathbb{C})=\left\{ F\colon \mathbb{C}\rightarrow \mathbb{C}%
\ \text{entire}:\Vert F\Vert _{\mathcal{F}^{2}(\mathbb{C})}^{2}=\int_{%
\mathbb{C}}|F(z)|^{2}e^{-\pi |z|^{2}}\,\mathrm{d}A(z)<\infty \right\} ,
\end{equation*}%
where $\mathrm{d}A(z)=dxdy$, for $z=x+iy$.\ This is a Hilbert space with
inner product 
\begin{equation*}
\left\langle F,G\right\rangle _{\mathcal{F}^{2}(\mathbb{C})}=\int_{\mathbb{C}%
}F(z)\overline{G(z)}e^{-\pi |z|^{2}}\,\mathrm{d}A(z),
\end{equation*}%
and reproducing kernel 
\begin{equation*}
K(z,w)=e^{\pi \overline{z}w}\text{.}
\end{equation*}%
The functions 
\begin{equation}
e_{n}(z)=\left( \frac{\pi ^{n}}{n!}\right) ^{1/2}z^{n},\qquad n=0,1,2,\ldots
,  \label{eq:monomial-basis}
\end{equation}%
form an orthonormal basis of $\mathcal{F}^{2}(\mathbb{C})$. For a measurable
set $\Omega \subset \mathbb{C}$ of finite area, let 
\begin{equation*}
T_{\Omega }F=P_{\mathcal{F}^{2}(\mathbb{C})}(\mathbf{1}_{\Omega }F)
\end{equation*}%
be the Toeplitz operator with symbol $\mathbf{1}_{\Omega }$, where $P_{%
\mathcal{F}^{2}(\mathbb{C})}$ is the orthogonal projection onto $\mathcal{F}%
^{2}(\mathbb{C})$. $T_{\Omega }$ can also be written as a Daubechies
time-frequency localization operator with Gaussian window \cite%
{Daubechies1988,Seip0,AbreDorf,AbrSpeck,NicolaTilli2022,GuerraGomesTilliRamos}%
. The quadratic form of $T_{\Omega }$ is: 
\begin{equation}
\left\langle T_{\Omega }F,F\right\rangle _{\mathcal{F}^{2}(\mathbb{C}%
)}=\int_{\Omega }|F(z)|^{2}e^{-\pi |z|^{2}}\,\mathrm{d}A(z).
\label{eq:quadratic-form-intro}
\end{equation}%
Since $\func{Tr}(T_{\Omega })=|\Omega |$, $T_{\Omega }$ is trace class
whenever $|\Omega |<\infty $. Moreover, it follows from (\ref%
{eq:quadratic-form-intro}) that $0\leq T_{\Omega }\leq I$. The eigenvalue $%
\lambda _{k}(\Omega )$\ associated to the eigenfunction $F_{k}^{\Omega }$\
can be written as%
\begin{equation*}
\lambda _{k}(\Omega )=\left\langle T_{\Omega }F_{k}^{\Omega },F_{k}^{\Omega
}\right\rangle _{\mathcal{F}^{2}(\mathbb{C})}=\int_{\Omega }|F_{k}^{\Omega
}(z)|^{2}e^{-\pi |z|^{2}}\,\mathrm{d}A(z)\text{.}
\end{equation*}%
The eigenvalues of $T_{\Omega }$ are arranged in nonincreasing order: 
\begin{equation*}
\lambda _{1}(\Omega )\geq \lambda _{2}(\Omega )\geq \cdots \geq 0\text{.}
\end{equation*}%
The largest eigenvalue provides the operator norm of $T_{\Omega }$ and
equals the optimal concentration of \emph{the Husimi function} of  $F\in 
\mathcal{F}^{2}(\mathbb{C})$ (known in time-frequency analysis as the \emph{%
spectrogram}). Nicola and Tilli proved that a disc maximizes $\lambda
_{1}(\Omega )$ among all measurable sets of a prescribed area \cite%
{NicolaTilli2022}. This is \emph{the Faber--Krahn inequality for the
short-time Fourier transform}, previously conjectured by the author and
Speckbacher \cite{AbrSpeck}. More recently, Nicola, Riccardi, and Tilli
proved the existence of a maximizing measurable set (without the radial
assumption used in the previous paper \cite[Theorem 1.3, Proposition 1.4]%
{NicolaRiccardiTilli2025}) such that 
\begin{equation}
\sup_{\lvert \Omega \rvert =s}\sum_{k=1}^{N}\lambda _{k}(\Omega )
\label{supKay}
\end{equation}%
is attained. The identification of that maximizer with a disc was left as a
conjecture in \cite[Page 2, after Theorem 1]{NicolaRiccardiTilli2026} (and
previously in \cite[Page 6]{NicolaRiccardiTilli2025}), and quite a
remarkable one: at a first glance one would not suspect that such a
statement could be true, since for individual eigenvalues with $k>1$ the
isoperimetric inequality fails. Indeed as shown in \cite[Theorem 1.3,
Proposition 1.4]{NicolaRiccardiTilli2025}, it already fails for $\lambda
_{2}(\Omega )$: among radial sets, the annulus is a maximizer, and even this
optimal property fails if the radial assumption is dropped. Thus, for the
conjecture to be true, a fine averaging of energy must be present, so that
each new eigenfunction added to the sum sharply fills in the energy
available, resulting in averages that are optimized by discs, independently
of the degree of the sum. In other words, this represents a \emph{completely
sharp version of the mechanism driving the accumulation of spectrograms} 
\cite{AGrRom,APR}. For radial sets, the Toeplitz operator is diagonal in the
monomial basis and the problem can be reduced to one-dimensional
inequalities. Nicola, Riccardi, and Tilli proved that discs maximize all
finite partial sums within this radial class \cite[Theorem 1.5]%
{NicolaRiccardiTilli2025}.

The problem cannot be handled by the modern tools developed after \cite%
{NicolaTilli2022} (see \cite%
{Kulikov,Frank,RamTilli,GuerraGomesTilliRamos,ProcLMS,CIMP, Ramos2026}),
that superseded the methods of \cite%
{Daubechies1988,dapa88,Seip0,AbreDorf,AbrSpeck,Galbis} (and certainly not by
these classical tools, that fail to deliver sharp bounds once radial
assumptions are relaxed). Indeed, already in \cite[page 3]%
{NicolaRiccardiTilli2025}, the authors raised the challenge (the brackets (%
\textbf{P. I})$_{\mathbf{I=1,2,3}}$ were inserted by the author):%
\begin{equation*}
\end{equation*}

`\emph{One reason of interest is that the machinery developed in \cite%
{NicolaTilli2022} for Faber--Krahn (though adaptable to contexts other than
the Fock space, see, e.g., \cite%
{Frank,FrankNicolaTilli,Kulikov,ortega,RamTilli}) is not suitable to attack
the aforementioned problems, since the differential inequality obtained in 
\cite{NicolaTilli2022} is no longer available: therefore, new ideas and
tools must be introduced. Moreover, contrary to (1.3), for the maximization
of }$\lambda _{2}(\Omega )$\emph{, circular symmetry breaks down, in that
the optimal set }$\Omega $\emph{\ (if one exists) is no longer rotational
invariant around any point. this symmetry breaking suggests that it might be
extremely difficult to characterize the sets }$\Omega $\emph{\ of prescribed
area \textsl{that maximize} }$\lambda _{2}(\Omega )$ (\textbf{P. 1})\emph{\
(or \textsl{the sum of the first k eigenvalues }}$\emph{(}$\textbf{P. 2}$%
\emph{)}$\emph{, or \textsl{the Schatten norms} }$\emph{(}$\textbf{P. 3}$%
\emph{)}$ \emph{etc.) without additional constraints, and these questions
remain (in their full generality) three challenging open problems.'}%
\begin{equation*}
\end{equation*}

We offer a solution to (\textbf{P. 2}) only assuming that $\left\vert \Omega
\right\vert <\infty $, a solution to (\textbf{P. 3}) and a bonus statement
regarding quasi-Schatten norms (these last two results follow from (\textbf{%
P. 2}) combined with the methods developed in \cite[page 3]%
{NicolaRiccardiTilli2025}). Regarding (\textbf{P. 1)}, we have no news
(shape maximization of individual eigenvalues for $k>1$ seems to offer
stubborn resistance to any approach). The main technical novelty is the 
incorporation of ideas from quantum information theory, which provides a
completely new perspective on how to handle shape optimization problems of
this nature. 

The main result of this paper is Theorem 1 below, showing that (\ref{supKay}%
) is maximized by discs of measure $s=\left\vert \Omega \right\vert $,
effectively confirming Nicola-Riccardi-Tilli conjecture from \cite%
{NicolaRiccardiTilli2026}. Moreover, combining the result with some of the
methods of \cite[page 3]{NicolaRiccardiTilli2025}, we will show that all
Schatten norms for $1<p<\infty $ are optimized by discs (a result that was
only previously known for $p=2$ \cite{NicolaRiccardi}; the question of
maximizing Schatten norms was raised in \cite[(1.5)]{NicolaRiccardiTilli2025}%
). We emphasize that, in \cite{NicolaRiccardiTilli2026}\ the authors refer
to (\ref{supKay}) as \emph{the Ky Fan norm.} Serendipitously enough, the
opening move of our argument, relating the Ky Fan norm of $T_{\Omega }$ (in
Petz's notation \cite[Theorem 11.13]{Petz}), 
\begin{equation}
\sigma _{N}\left( T_{\Omega }\right) =\sum_{k=1}^{N}\lambda _{k}(\Omega )%
\text{,}  \label{KF}
\end{equation}%
with the maximal Husimi concentration among rank $N$ projection operators,
relies on Ky Fan's maximum principle \cite[Theorem 1, formula (4)]{Fan1951}.
We now state the main result. For $s>0$, write $D_{s}$ for the centered disc
of area $s$, 
\begin{equation*}
D_{s}=\{z\in \mathbb{C}:\pi |z|^{2}<s\}\text{.}
\end{equation*}

\medskip \noindent \textbf{Theorem 1 (Isoperimetric inequality for the Ky
Fan norm).} \emph{Let }$\Omega \subset \mathbb{C}$\emph{\ be measurable with }$%
s=|\Omega |<\infty $\emph{, and let }$N\geq 1$\emph{. Then }%
\begin{equation}
\sum_{j=1}^{N}\lambda _{j}(\Omega )\leq \sum_{j=1}^{N}\lambda
_{j}(D_{s})=N-e^{-|\Omega |}\sum_{k=0}^{N-1}(N-k)\frac{|\Omega |^{k}}{k!}.
\label{eq:main}
\end{equation}%
\emph{The bound is attained by the centered disc of area }$s$\emph{. By Weyl
translation invariance, it is also attained by every disc of area }$s$\emph{.%
} \medskip

The only assumption on $\Omega $ used in Theorem 1 is finite area. No
smoothness, boundedness, or connectedness assumption on $\Omega $ are
required. For radial sets Theorem 1 was proved in \cite[Theorem 1.5]%
{NicolaRiccardiTilli2025}. The above inequality has far-reaching
consequences. When $N=1$, this gives%
\begin{equation*}
\lambda _{1}(\Omega )\leq 1-e^{-|\Omega |}
\end{equation*}%
and recovers the result of Nicola and Tilli \cite{NicolaTilli2022}. For $N=2$%
, the sharp bound is $\lambda _{1}(\Omega )+\lambda _{2}(\Omega )\leq
2-(2+|\Omega |)e^{-|\Omega |}$. 

\begin{remark}
Nicola, Riccardi and Tilli \cite{NicolaRiccardiTilli2026} proved the
existence of an optimal set and conjectured that discs are optimal. The
uniqueness of optimal sets, up to translations and null sets, was left open.
Our argument does not address this question.
\end{remark}

As a significant application of Theorem 1, we obtain the isoperimetric
Schatten inequalities stated in Corollary 1 below. We proceed with the
definitions of Schatten norms and quasi norms. For $0<p<\infty $, the
Schatten $p$-sum of a positive compact operator $T$ is 
\begin{equation*}
\Vert T\Vert _{S^{p}}=\left( \sum_{j=1}^{\infty }\lambda _{j}(T)^{p}\right)
^{1/p}.
\end{equation*}%
For $p\geq 1$ $\Vert T\Vert _{S^{p}}$ is a norm. If $T_{\Omega }\in S^{p}$
for $0<p<1$, $\Vert T\Vert _{S^{p}}$ is a quasinorm. We also write 
\begin{equation*}
\Vert T\Vert _{S^{\infty }}=\lambda _{1}(T)=\Vert T\Vert _{\mathrm{op}}.
\end{equation*}%
The following result extends \cite[Corollary~1.8]{NicolaRiccardiTilli2025}
from circularly symmetric sets to arbitrary measurable sets, and extends 
\cite{NicolaRiccardi} from $p=2$ to general $0<p<\infty .$ It bounds all
possible ranges of Schatten norms. The proof of Corollary 1 follows
combining Theorem 1 with the infinite-dimensional Karamata principle from 
\cite[Corollary 3.4, Proposition~3.5]{NicolaRiccardiTilli2025}.

\medskip \noindent \textbf{Corollary 1 (Isoperimetric Schatten inequalities).%
} \emph{Let }$\Omega \subset \mathbb{C}$\emph{\ be measurable with }$|\Omega
|=s<\infty $\emph{. Then: }%
\begin{equation*}
\begin{array}{ll}
\Vert T_{\Omega }\Vert _{S^{p}}\leq \Vert T_{D_{s}}\Vert
_{S^{p}}=\left[\displaystyle\sum_{n=0}^{\infty }\left(
1-e^{-s}\sum_{k=0}^{n}\frac{s^{k}}{k!}\right) ^{p}\right]^{1/p}
& 1<p<\infty , \\[3mm]
\Vert T_{\Omega }\Vert _{S^{\infty}}\leq
\Vert T_{D_s}\Vert _{S^{\infty}}=1-e^{-s},
& p=\infty, \\[2mm]
\Vert T_{\Omega }\Vert _{S^{1}}=\Vert T_{D_{s}}\Vert _{S^{1}}=|\Omega |, & 
p=1, \\[2mm] 
\func{Tr}(T_{\Omega }^{p})\geq \func{Tr}(T_{D_{s}}^{p})=\sum_{n=0}^{\infty
}\left( 1-e^{-s}\sum_{k=0}^{n}\frac{s^{k}}{k!}\right) ^{p}
& 0<p<1.%
\end{array}%
\end{equation*}%
\emph{The bounds are attained by every disc of area }$s$\emph{. The last
inequality holds for the quasi-Schatten norms if }$T_{\Omega }\in S^{p}$%
\emph{, in which case }%
\begin{equation*}
\Vert T_{\Omega }\Vert _{S^{p}}\geq \Vert T_{D_{s}}\Vert _{S^{p}},\qquad
0<p<1.
\end{equation*}

Our proof of Theorem 1 takes a surprisingly simple and direct approach. The
key observation is the following: while \emph{the largest eigenvalue
maximizes the concentration of Husimi functions} $Q_{F}$, where $F$ runs
over $\mathcal{F}^{2}(\mathbb{C})$ (or over all rank $1$ operators $%
|F\rangle \langle F|$ in $\mathcal{F}^{2}(\mathbb{C})$), \emph{the sum of
the first }$N$\emph{\ eigenvalues maximizes the concentration of Husimi
functions }$Q_{\rho _{P}}$, when $\rho _{P}=P/N$ and $P$\ runs over all rank-%
$N$ orthogonal projections in $\mathcal{F}^{2}(\mathbb{C})$. This last
property will be established using Ky Fan's maximum principle (\cite[Theorem
1, formula (4)]{Fan1951},\cite[Theorem 11.13]{Petz}). This allows us to use
results from quantum information theory, namely De Palma's majorization
theorem for the Husimi transform of density operators \cite[Theorem~5]%
{DePalma2018}, therefore resorting to the theory of Fock rearrangements of
density operators put forward by De Palma, Trevisan and Giovannetti \cite%
{PslmsTrevissanGiovannetti}. An adaptation to $\mathbb{C}$ of the bathtub
principle from \cite{AbrSpeck} provides the third tool used in the argument. 

Despite the short proof, Theorem 1 is a profound result. Indeed, in \cite%
{DePalma2018}, De Palma uses several deep results to obtain the majorization
theorem. In particular, the argument relies on operator versions of the Lieb
inequalities of \cite{Lieb1978,Lieb1990} and a deep result of \cite[Theorem
V.3]{PslmsTrevissanGiovannetti}, which states that, among all positive
operators with prescribed eigenvalues, placing those eigenvalues in
decreasing order along the Fock basis maximizes every finite partial sum of
the output eigenvalues.

\subsection{Previous related work}

This result is one more contribution to a series of developments related to
the interval conjecture for the Paley-Wiener space, recently confirmed
in \cite{AS2026}. This conjecture was formulated in the influential paper of
Donoho and Stark \cite[Conjecture 1]{DonohoStark}, which, together with the $%
L_{1}$ large sieve of Donoho and Logan \cite{DonohoLogan}, laid the
foundations for the mathematical understanding of the \emph{%
sparsity/concentration trade-off} characteristic of signal recovering
methods.

The concentration part of the sparsity/concentration trade-off has attracted
significant attention recently, following analogies with the Laplacian
problem on a domain. The interval concentration conjecture of Donoho and
Stark was reformulated in \cite[Conjecture 1]{AbrSpeck} for
Daubechies localization operators \cite{Daubechies1988} (where the phase
space can be identified with the Fock space of entire functions \cite{Seip0}%
)\ and for Daubechies-Paul localization operators \cite{dapa88} in \cite%
{AbreDorf} (where the phase space can be identified with the Bergman space
of analytic functions \cite{Seip0}). After some preliminary advances \cite%
{Galbis}, a breakthrough in the Gaussian Short Time Fourier Transform/Fock
case was achieved in \cite{NicolaTilli2022}. The methods of \cite%
{NicolaTilli2022} proved to be effective for spaces of analytic functions
with constant curvature, independently of the geometry involved, leading to
a solution of the problem for Cauchy Wavelets/Bergman space \cite{RamTilli},
Hardy space \cite{ProcLMS} and more general coherent state phase spaces \cite%
{Frank} (including analytic polynomials on the sphere). The rigidity of
concentration on discs exhibited by the inverse \cite{AbreDorf} and optimal
concentration problems \cite{NicolaTilli2022} has been recently explored
from different perspectives, including free boundary problems \cite%
{Ramos2026} and sharp geometric analysis stability \cite%
{GuerraGomesTilliRamos}, with extensions to the hyperbolic \cite{ProcLMS}
and spherical settings \cite{CIMP}. Other related advances stem from the
observation that \cite{NicolaTilli2022} offers an alternative proof of
Lieb's inequality for the Fock space \cite{NicolaTilli2023}. This subtlety
was previously observed by Kulikov \cite{Kulikov}, who adapted the methods
to prove the Lieb inequality for the Cauchy wavelets/Bergman spaces.

\subsection{Outline of the paper}

The paper is organized as follows. Section~\ref{sec:fock} introduces the
Husimi function of an operator and derives the Ky Fan representation.
Section~\ref{sec:majorization} states De Palma's theorem in the
normalization used here. Section~\ref{sec:proof} recalls the bathtub
principle, proves Theorem 1, and evaluates the sharp constant. Finally,
Section~\ref{sec:schatten} proves Corollary 1.

\section{The Ky Fan representation for the Husimi function}

\label{sec:fock}

\subsection{A Rayleigh quotient for the Ky Fan norm}

It will be convenient to use the weighted reproducing kernel, 
\begin{equation}
k_{z}(w)=e^{-\pi |z|^{2}/2}K(z,w)=e^{\pi \overline{z}w-\pi |z|^{2}/2}\text{,}
\label{eq:coherent}
\end{equation}%
which is normalized such that $\Vert k_{z}\Vert _{\mathcal{F}^{2}(\mathbb{C}%
)}=1$. The reproducing kernel formula then gives 
\begin{equation}
|\left\langle F,k_{z}\right\rangle _{\mathcal{F}^{2}(\mathbb{C}%
)}|^{2}=|F(z)|^{2}e^{-\pi |z|^{2}}\text{.}  \label{eq:normalized-reproducing}
\end{equation}%
A density operator is a positive trace-class operator $\rho $ with $\func{Tr}%
(\rho )=1$ \cite[Section 2.1]{Petz}. Its Husimi function is defined as 
\begin{equation}
Q_{\rho }(z)=\left\langle k_{z},\rho k_{z}\right\rangle _{\mathcal{F}^{2}(%
\mathbb{C})}\text{.}  \label{eq:Husimi}
\end{equation}%
The Husimi function is nonnegative and satisfies 
\begin{equation}
0\leq Q_{\rho }(z)\leq 1,\qquad \int_{\mathbb{C}}Q_{\rho }(z)\,\mathrm{d}%
A(z)=1\text{.}  \label{eq:Husimi-normalization}
\end{equation}%
The Husimi function may therefore be viewed as an ordinary probability
density on phase space. Now consider the density operator 
\begin{equation}
\rho _{P}=\frac{1}{N}P\text{,}  \label{densityOperator}
\end{equation}%
where $P$ is an orthogonal projection of rank $N$ in $\mathcal{F}^{2}(%
\mathbb{C})$. Its Husimi function is 
\begin{equation}
Q_{\rho _{P}}(z)=\frac{1}{N}\left\langle k_{z},Pk_{z}\right\rangle _{%
\mathcal{F}^{2}(\mathbb{C})}\text{.}  \label{eq:qP}
\end{equation}%
If $F_{1},\ldots ,F_{N}$ is an orthonormal basis of the range of $P$, then
\begin{equation*}
Pk_z=\sum_{j=1}^{N}\langle k_z,F_j\rangle_{\mathcal{F}^{2}(\mathbb{C})}F_j,
\end{equation*}
and therefore
\begin{equation}
Q_{\rho _{P}}(z)=\frac{1}{N}\sum_{j=1}^{N}|\left\langle
k_{z},F_{j}\right\rangle _{\mathcal{F}^{2}(\mathbb{C})}|^{2}=\frac{1}{N}%
e^{-\pi |z|^{2}}\sum_{j=1}^{N}|F_{j}(z)|^{2}\text{.}  \label{eq:qP-basis}
\end{equation}%
This shows, in particular, that $Q_{\rho_P}$ depends only on the range of
$P$, not on the chosen orthonormal basis. Since $P$ is an orthogonal
projection,
\begin{equation*}
Q_{\rho_P}(z)=\frac{1}{N}\langle k_z,Pk_z\rangle
=\frac{1}{N}\Vert Pk_z\Vert^2\leq\frac{1}{N}.
\end{equation*}
Thus,
\begin{equation*}
0\leq Q_{\rho _{P}}\leq \frac{1}{N}\text{.}
\end{equation*}%
Integrating (\ref{eq:qP-basis}) gives 
\begin{equation}
\int_{\mathbb{C}}Q_{\rho _{P}}\,\mathrm{d}A(z)=\frac{1}{N}%
\sum_{j=1}^{N}\int_{\mathbb{C}}|F_{j}(z)|^{2}e^{-\pi |z|^{2}}\mathrm{d}A(z)=1%
\text{.}  \label{eq:qP-mass}
\end{equation}

Then (\ref{eq:qP-basis}) combined with the next Proposition shows that $%
\frac{1}{N}\sum_{j=1}^{N}\lambda _{j}(\Omega )$ is exactly the Rayleigh
quotient of $Q_{\rho _{P}}$.

\medskip \noindent \textbf{Proposition 2.1 (Ky Fan representation).} For
every measurable $\Omega \subset \mathbb{C}$ of finite area, 
\begin{equation}
\frac{1}{N}\sum_{j=1}^{N}\lambda _{j}(\Omega )=\max_{\substack{ P=P^{\ast
}=P^{2}  \\ \func{rank}(P)=N}}\int_{\Omega }Q_{\rho _{P}}\,\mathrm{d}A(z)%
\text{,}  \label{eq:kyfan}
\end{equation}%
where the maximum is taken over all rank-$N$ orthogonal projections $P$ and $%
\rho _{P}=\frac{1}{N}P$.

\begin{proof}
The operator $T_{\Omega }$ is positive and compact, with eigenvalues in
nonincreasing order, 
\begin{equation*}
\lambda _{1}(\Omega )\geq \lambda _{2}(\Omega )\geq \cdots \geq 0\text{.}
\end{equation*}%
Take $m=1$, $A_{1}=T_{\Omega }$ and $U_{1}=I$\ in Ky Fan's maximum principle 
\cite[Theorem 1, formula (4)]{Fan1951} (the precise form of the result we
use is stated in \cite[Theorem 11.13]{Petz} without proof). Since $T_{\Omega
}$ is positive and compact, its eigenvalues coincide with the singular
values and 
\begin{equation*}
\sum_{j=1}^{N}\left\langle T_{\Omega }F_{j},F_{j}\right\rangle _{\mathcal{F}%
^{2}(\mathbb{C})}\leq \sum_{j=1}^{N}\lambda _{j}(\Omega )\text{,}
\end{equation*}%
for every family $\{F_{1},...,F_{N}\}$ of orthonormal functions in the Fock
space. Equality is attained when $F_{1},...,F_{N}$ are eigenfunctions
corresponding to the first $N$ eigenvalues. Thus,%
\begin{equation*}
\sum_{j=1}^{N}\lambda _{j}(\Omega )=\max_{\substack{ F_{1},\ldots ,F_{N}\in 
\mathcal{F}^{2}(\mathbb{C}) \\ \left\langle F_{j},F_{k}\right\rangle _{%
\mathcal{F}^{2}(\mathbb{C})}=\delta _{jk}}}\sum_{j=1}^{N}\left\langle
T_{\Omega }F_{j},F_{j}\right\rangle _{\mathcal{F}^{2}(\mathbb{C})}\text{.}
\end{equation*}%
Now, let $P$ be the orthogonal projection onto $\operatorname{span}\{F_{1},\ldots ,F_{N}\}
$.\ For this projection,%
\begin{equation*}
\func{Tr}(PT_{\Omega })=\sum_{j=1}^{N}\left\langle T_{\Omega
}F_{j},F_{j}\right\rangle _{\mathcal{F}^{2}(\mathbb{C})}\text{.}
\end{equation*}%
Therefore, maximizing over orthonormal families is equivalent to maximizing
over rank-$N$ orthogonal projections:%
\begin{equation*}
\sum_{j=1}^{N}\lambda _{j}(\Omega )=\max_{\substack{ P=P^{\ast }=P^{2} \\ 
\func{rank}(P)=N}}\func{Tr}(PT_{\Omega })\text{.}
\end{equation*}%
Using (\ref{eq:quadratic-form-intro}) and (\ref{eq:qP-basis}), we have, for
the density operator (\ref{densityOperator}), 
\begin{equation*}
\func{Tr}(\rho _{P}T_{\Omega })=\frac{1}{N}\sum_{j=1}^{N}\left\langle T_{\Omega
}F_{j},F_{j}\right\rangle _{\mathcal{F}^{2}(\mathbb{C})}=\frac{1}{N}%
\int_{\Omega }e^{-\pi |z|^{2}}\sum_{j=1}^{N}|F_{j}(z)|^{2}\,\mathrm{d}%
A(z)=\int_{\Omega }Q_{\rho _{P}}(z)\mathrm{d}A(z)\text{.}
\end{equation*}%
Conversely, every rank-$N$ projection arises from an orthonormal basis of
its range. This proves (\ref{eq:kyfan}).
\end{proof}

\section{De Palma's Husimi majorization}

\label{sec:majorization}

In \cite{PslmsTrevissanGiovannetti}, De Palma, Trevisan, and Giovannetti
used Fock rearrangements to show that the output generated by any quantum
state is majorized by the output generated by the state with the same
spectrum diagonal in the Fock basis and with decreasing eigenvalues (see
also \cite[Theorem 1]{DePalma2018}). This was a crucial result in the proof
of our main tool, \cite[Theorem 5]{DePalma2018} (which also depends on the
Lieb inequality for density operators \cite[Theorems 8,9]{DePalma2018}). De
Palma's Husimi majorization theorem is presented in this section, after
defining Fock rearrangements.

\begin{definition}
Consider the density operator%
\begin{equation*}
\rho =\sum_{n=1}^{\infty }p_{n}|F_{n}\rangle \langle F_{n}|,\quad p_{1}\geq
p_{2}\geq \ldots \geq 0,\quad \langle F_{m}|F_{n}\rangle =\delta _{mn}\text{.%
}
\end{equation*}%
Its \emph{Fock }(in the terminology of \cite{PslmsTrevissanGiovannetti}) or%
\emph{\ passive rearrangement} (in the terminology of \cite[Theorem 5]%
{DePalma2018}) is defined as the operator with the same spectrum as $\rho $
, and given as 
\begin{equation}
\rho ^{\downarrow }=\sum_{n=1}^{\infty }p_{n}|e_{n-1}\rangle \langle e_{n-1}|%
\text{,}  \label{eq:passive}
\end{equation}%
where%
\begin{equation*}
p_{1}\geq p_{2}\geq \cdots \geq 0
\end{equation*}%
are the eigenvalues of $\rho $, repeated according to multiplicity. We say
that $\hat{\rho}$ is \emph{passive} if $\hat{\rho}=\hat{\rho}^{\downarrow }$%
, i.e. if $\hat{\rho}$ is diagonal in the Fock basis with eigenvalues
decreasing as the energy increases.
\end{definition}

Our key instrument will be the following theorem of De Palma \cite[Theorem 5]%
{DePalma2018}. We state only the form required for the proof of Theorem 1.

\medskip \noindent \textbf{Theorem 3.1 (De Palma's Husimi majorization
theorem).} \emph{Let }$\rho $\emph{\ be a density operator on }$\mathcal{F}^{2}(\mathbb{C})$%
\emph{\ and let }$\rho ^{\downarrow }$\emph{\ be its passive rearrangement.
If }$\Phi \in C^{1}([0,1])$\emph{\ is convex and }$\Phi (0)=0$\emph{, then }%
\begin{equation}
\int_{\mathbb{C}}\Phi (Q_{\rho }(z))\,\mathrm{d}A(z)\leq \int_{\mathbb{C}%
}\Phi (Q_{\rho ^{\downarrow }}(z))\,\mathrm{d}A(z).  \label{eq:depalma}
\end{equation}

In \cite[Theorem~5]{DePalma2018} the coherent state is indexed by $\alpha
\in \mathbb{C}$ and the phase-space measure is$\,\mathrm{d}^{2}\alpha /\pi $%
. The change of variables $\alpha =\sqrt{\pi }\,z$ gives exactly (\ref%
{eq:depalma}). De Palma's theorem states that, among all states with a fixed
spectrum, the one with eigenvalues ordered decreasingly has the most
concentrated Husimi function in the sense of (\ref{eq:depalma}). Before
proving the main result we need the following Proposition, which is a simple
application of Theorem 3.1 to a rank-$N$ projection.

\medskip \noindent \textbf{Proposition 3.2.} \emph{Let }$P$\emph{\ be any
orthogonal projection of rank }$N$\emph{\ and }$\rho =P/N$\emph{\ its
associated density operator. Then, for every }$\tau \geq 0$\emph{, }%
\begin{equation}
\int_{\mathbb{C}}(Q_{\rho }(z)-\frac{\tau }{N})_{+}\,\mathrm{d}A(z)\leq
\int_{\mathbb{C}}(Q_{\rho ^{\downarrow }}(z)-\frac{\tau }{N})_{+}\,\mathrm{d}%
A(z)\text{.}  \label{eq:hinge-majorization}
\end{equation}

\begin{proof}
For $0<\tau <1$, consider the convex function 
\begin{equation*}
\Phi (x)=\left( x-\frac{\tau }{N}\right) _{+}.
\end{equation*}%
Since $\Phi $ is not $C^{1}$ at $\tau /N$, a standard approximation is
required for the application of Theorem 3.1. For $a=\tau /N$ and $%
\varepsilon >0$, define 
\begin{equation*}
\Phi _{\varepsilon }(x)=%
\begin{cases}
0, & (x-a)\leq 0, \\[2mm] 
\dfrac{(x-a)^{2}}{2\varepsilon }, & 0<(x-a)<\varepsilon , \\[3mm] 
(x-a)-\dfrac{\varepsilon }{2}, & (x-a)\geq \varepsilon .%
\end{cases}%
\end{equation*}%
For every $\varepsilon >0$, $\Phi _{\varepsilon }\in C^{1}([0,1])$.
Moreover, $\Phi _{\varepsilon }$ is convex, with $\Phi _{\varepsilon }(0)=0$
and $\lim\limits_{\epsilon \rightarrow 0}\Phi _{\varepsilon }(x)=(x-a)_{+}$
pointwise. Since $0\leq \Phi _{\varepsilon }(x)\leq x$, (\ref%
{eq:Husimi-normalization}) shows that the convergence is dominated. Applying
Theorem 3.1 and passing to the limit, gives 
\begin{equation*}
\int_{\mathbb{C}}(Q_{\rho }(z)-\frac{\tau }{N})_{+}\mathrm{d}A(z)\leq \int_{%
\mathbb{C}}(Q_{\rho ^{\downarrow }}(z)-\frac{\tau }{N})_{+}\mathrm{d}A(z).
\end{equation*}%
This is (\ref{eq:hinge-majorization}). When $\tau =0$, both sides equal $1$
by (\ref{eq:qP-mass}). When $\tau \geq 1$, both sides vanish because $0\leq
Q_{\rho },Q_{\rho ^{\downarrow }}\leq 1/N\leq \frac{\tau }{N}$.
\end{proof}

\section{Proof of Theorem 1}

\label{sec:proof}

We will use the following form of the bathtub principle, which is an
adaptation of \cite[Lemma 2.3]{AS2026} to $\mathbb{C}$. Similar results are
proved in \cite{LiebLoss2001}.

\noindent \textbf{Lemma (Bathtub principle).} Let $v\colon \mathbb{C}%
\rightarrow \lbrack 0,\infty )$ be integrable and let $s>0$. Choose $\tau
\geq 0$ such that 
\begin{equation*}
|\{v>\tau \}|\leq s\leq |\{v\geq \tau \}|,
\end{equation*}%
and choose a measurable set $F$ satisfying 
\begin{equation*}
\{v>\tau \}\subseteq F\subseteq \{v\geq \tau \},\qquad |F|=s.
\end{equation*}%
Then, for every measurable set $A\subseteq \mathbb{C}$ with $|A|=s$, 
\begin{equation*}
\int_{A}v\,\,\mathrm{d}A\leq \int_{F}v\,\,\mathrm{d}A.
\end{equation*}

\textbf{Proof of Theorem 1}

The operator $\rho =P/N$ is a density operator whose spectrum is $%
\underbrace{\{\frac{1}{N},\ldots ,\frac{1}{N}}_{N\ \mathrm{times}%
},0,0,\ldots .\}$. Its passive rearrangement is therefore 
\begin{equation*}
\rho ^{\downarrow }=\rho _{N}=\frac{1}{N}\sum_{n=0}^{N-1}|e_{n}\rangle \langle
e_{n}|.
\end{equation*}%
We can compute $Q_{\rho _{N}}$ explicitly. From (\ref{eq:monomial-basis})
and (\ref{eq:normalized-reproducing}), 
\begin{equation*}
Q_{\rho _{N}}(z)=\frac{1}{N}\sum_{n=0}^{N-1}|\left\langle
e_{n},k_{z}\right\rangle _{\mathcal{F}^{2}(\mathbb{C})}|^{2}=\frac{1}{N}%
e^{-\pi |z|^{2}}\sum_{n=0}^{N-1}\frac{(\pi |z|^{2})^{n}}{n!}\text{.}
\end{equation*}%
Setting $t=\pi |z|^{2}$ and differentiating, 
\begin{equation}
\frac{\,\mathrm{d}}{\,\mathrm{d}t}\left( \frac{1}{N}e^{-t}\sum_{n=0}^{N-1}%
\frac{t^{n}}{n!}\right) =-e^{-t}\frac{t^{N-1}}{N!}.  \label{eq:qN-derivative}
\end{equation}%
Thus, $Q_{\rho _{N}}(z)$ is radial and strictly decreasing away from the
origin. It follows from the bathtub principle that the set of area $s$
capturing the largest amount of $Q_{\rho _{N}}$ is the centered disc $D_{s}$
of area $s=|\Omega |$. Writing 
\begin{equation*}
r_{s}=\sqrt{\frac{s}{\pi }}\qquad \text{and \ \ \ \ }\tau _{s}=Q_{\rho
_{N}}(r_{s})\text{,}
\end{equation*}%
the \emph{open} disc $D_{s}$ can be written as the superlevel set 
\begin{equation*}
D_{s}=\{z\in \mathbb{C}:Q_{\rho _{N}}(z)>\tau _{s}\}
\end{equation*}%
up to its boundary. This implies that%
\begin{equation*}
\left( Q_{\rho _{N}}(z)-\tau _{s}\right) _{+}=\left\{ 
\begin{array}{cc}
Q_{\rho _{N}}(z)-\tau _{s}\text{, } & z\in D_{s} \\ 
0 & z\notin D_{s}%
\end{array}%
\right. 
\end{equation*}%
almost everywhere. Thus,%
\begin{equation}
\int_{\mathbb{C}}\left( Q_{\rho _{N}}(z)-\tau _{s}\right) _{+}\,\mathrm{d}%
A(z)=\int_{D_{s}}\left( Q_{\rho _{N}}(z)-\tau _{s}\right) \,\mathrm{d}A(z)
\label{Identity}
\end{equation}%
Since $|\Omega |=s$, this gives 
\begin{eqnarray*}
\int_{\Omega }Q_{\rho }(z)\,\mathrm{d}A(z) &=&s\tau _{s}+\int_{\Omega
}\left( Q_{\rho }(z)-\tau _{s}\right) \,\mathrm{d}A(z) \\
&\leq &s\tau _{s}+\int_{\mathbb{C}}\left( Q_{\rho }(z)-\tau _{s}\right)
_{+}\,\mathrm{d}A(z) \\
&\leq &s\tau _{s}+\int_{\mathbb{C}}(Q_{\rho _{N}}(z)-\tau _{s})_{+}\,\mathrm{%
d}A(z)\text{.} \\
&=&s\tau _{s}+\int_{D_{s}}(Q_{\rho _{N}}(z)-\tau _{s})_{+}\,\mathrm{d}A(z) \\
&=&\int_{D_{s}}Q_{\rho _{N}}(z)\,\mathrm{d}A(z)\text{,}
\end{eqnarray*}%
where the inequality in the third line follows from Proposition 3.2 with $%
\tau =N\tau _{s}$ and the last identity uses $|D_{s}|=s$.

Finally, multiplying by $N$, taking the supremum over $P$ and rewriting the
left hand side of the inequality using the Ky Fan representation from
Proposition 2.1, gives 
\begin{equation}
\sum_{j=1}^{N}\lambda _{j}(\Omega )\leq N\int_{D_{s}}Q_{\rho _{N}}(z)\,%
\mathrm{d}A(z).  \label{eq:main-before-evaluation}
\end{equation}%
It only remains to show that the right-hand side is exactly the sum of the
first $N$ eigenvalues for the disc. This is a routine calculation.
Rotational invariance shows that $\{e_{n}\}_{n\geq 0}$ is orthogonal in any
centered disc (for details see, e.g., the proof of \cite[Theorem 3.4.2 a)]%
{Grochenig2001}):%
\begin{equation}
\left\langle T_{D_{s}}e_{n},e_{m}\right\rangle _{\mathcal{F}^{2}(\mathbb{C}%
)}=\int_{D_{s}}e_{n}(z)\overline{e_{m}(z)}e^{-\pi |z|^{2}}\,\mathrm{d}%
A(z)=\delta _{n,m}\lambda _{n+1}(D_{s})\text{.}  \label{double}
\end{equation}%
Using polar coordinates and the substitution $t=\pi |z|^{2}$, (\ref{double})
implies that 
\begin{equation}
\lambda _{n+1}(D_{s}):=\left\langle T_{D_{s}}e_{n},e_{n}\right\rangle _{%
\mathcal{F}^{2}(\mathbb{C})}=\frac{1}{n!}\int_{0}^{s}e^{-t}t^{n}\,\mathrm{d}%
t=1-e^{-s}\sum_{k=0}^{n}\frac{s^{k}}{k!}  \label{discEigenvalues}
\end{equation}%
is the eigenvalue of $T_{D_{s}}$ associated with $e_{n}$ (this follows from
well known argument, see \cite[Proposition 1]{Seip0}, or \cite[Proposition
4.2, Proposition 5.1]{Svela}, where this property was extended to mixed
state localization operators). Thus, $T_{D_{s}}$ is diagonal in the basis $%
\{e_{n}\}_{n\geq 0}$ . Moreover, 
\begin{equation*}
\lambda _{n+2}(D_{s})-\lambda _{n+1}(D_{s})=-e^{-s}\frac{s^{n+1}}{(n+1)!}<0%
\text{,}
\end{equation*}%
showing that $\lambda _{1}(D_{s}),\ldots ,\lambda _{N}(D_{s})$ are the
first $N$ eigenvalues of $T_{D_{s}}$ arranged in decreasing order.
Consequently, (\ref{double}) implies:%
\begin{equation}
\int_{D_{s}}Q_{\rho _{N}}(z)\,\mathrm{d}A(z)=\frac{1}{N}\sum_{n=0}^{N-1}%
\int_{D_{s}}|e_{n}(z)|^{2}e^{-\pi |z|^{2}}\mathrm{d}A(z)=\frac{1}{N}%
\sum_{n=1}^{N}\lambda _{n}(D_{s})  \label{HusimieigengaluesDisc}
\end{equation}%
Finally, summing (\ref{discEigenvalues}), 
\begin{equation*}
\sum_{n=1}^{N}\lambda _{n}(D_{s})=\sum_{n=1}^{N}\left(
1-e^{-s}\sum_{k=0}^{n-1}\frac{s^{k}}{k!}\right) =N-e^{-\left\vert \Omega
\right\vert }\sum_{n=0}^{N-1}(N-n)\frac{\left\vert \Omega \right\vert ^{n}}{%
n!}\text{.}
\end{equation*}%
Combining this with (\ref{eq:main-before-evaluation}) and (\ref%
{HusimieigengaluesDisc}), gives%
\begin{equation*}
\sum_{j=1}^{N}\lambda _{j}(\Omega )\leq N\int_{D_{s}}Q_{\rho _{N}}(z)\,%
\mathrm{d}A(z)=\sum_{n=1}^{N}\lambda _{n}(D_{s})=N-e^{-\left\vert \Omega
\right\vert }\sum_{n=0}^{N-1}(N-n)\frac{\left\vert \Omega \right\vert ^{n}}{n!}%
\text{,}
\end{equation*}%
which is exactly (\ref{eq:main}). Invariance under Weyl translations gives
equality for every disc with area equal to $\left\vert \Omega \right\vert $.

\section{Proof of Corollary 1}
\label{sec:schatten}

We will use the concepts of weak and strong majorization \cite%
{HardyLittlewoodPolya,Petz}. Given two nonincreasing summable sequences
$x=(x_{j})_{j=1}^{\infty }$ and $y=(y_{j})_{j=1}^{\infty }$ with
nonnegative entries, $x$ is said to be \emph{weakly majorized} by $y$,
written $x\prec _{w}y$, if
\begin{equation}
\sum_{j=1}^{N}x_{j}\leq \sum_{j=1}^{N}y_{j}
\qquad\text{for every }N\geq1. \label{xsuccwy}
\end{equation}%
It is \emph{strongly majorized} by $y$, written $x\prec y$, if
$x\prec_w y$ and
\begin{equation*}
\sum_{j=1}^{\infty}x_j=\sum_{j=1}^{\infty}y_j.
\end{equation*}

In this context, Theorem 1 gives, for every $N<\infty $, 
\begin{equation*}
\sum_{j=1}^{N}\lambda _{j}(\Omega )\leq \sum_{j=1}^{N}\lambda _{j}(D_{s})%
\text{,}
\end{equation*}%
while, for $N=\infty $, 
\begin{equation*}
\sum_{j=1}^{\infty }\lambda _{j}(\Omega )=\func{Tr}(T_{\Omega })=|\Omega
|=s=\sum_{j=1}^{\infty }\lambda _{j}(D_{s})\text{.}
\end{equation*}%
Thus, we have the strong majorization relation 
\begin{equation}
\lambda (T_{\Omega })\prec \lambda (T_{D_{s}}),\qquad |\Omega |=s.
\label{eq:spectral-majorization}
\end{equation}

Now we apply to (\ref{eq:spectral-majorization}) the infinite-dimensional
Karamata principle from \cite[Proposition~3.5]{NicolaRiccardiTilli2025} to
yield the proof.

\begin{proof}
For $1<p<\infty $, apply \cite[Corollary 3.4]{NicolaRiccardiTilli2025} to
the convex function $\Phi (t)=t^{p}$. For $0<p<1$, apply the
infinite-dimensional Karamata inequality from \cite[Proposition~3.5]%
{NicolaRiccardiTilli2025} to the convex function $\Phi (t)=-t^{p}$. The case 
$p=1$ is $\func{Tr}(T_{\Omega })=|\Omega |$, and the case $p=\infty $ is the 
$N=1$ case of Theorem 1. The sharp constants are explicit. From (\ref%
{discEigenvalues}), for every $0<p<\infty $, 
\begin{equation}
\Vert T_{D_{s}}\Vert _{S^{p}}^{p}=\sum_{n=0}^{\infty }\left(
1-e^{-s}\sum_{k=0}^{n}\frac{s^{k}}{k!}\right) ^{p}.
\label{eq:sharp-schatten-constant}
\end{equation}%
The series converges for every $p>0$, because of the factorial decay of
the eigenvalues. At the endpoints, 
\begin{equation*}
\Vert T_{D_{s}}\Vert _{S^{1}}=s,\qquad \Vert T_{D_{s}}\Vert _{S^{\infty
}}=1-e^{-s}.
\end{equation*}
\end{proof}

\begin{remark}
These methods require analyticity due to the structure of maximizers in de
Palma's theorem. For general Toeplitz operators with translation invariant
kernels or for Daubechies operators with general windows, restricted two
sided inequalities were proved in \cite{APR} (not explicitly, but can be
derived from the last sections), with the goal of obtaining the sharp rate
of $L_{1}$ deviations of the accumulated spectrogram with a general window $g
$ (the setting considered in this paper follows from taking $g$ a Gaussian
function):%
\begin{equation*}
\Vert T_{\Omega }\Vert _{S^{2}}\leq \sum_{k=1}^{N_{\Omega }}\lambda
_{k}^{\Omega }\leq \frac{N_{\Omega }}{2}-\frac{\Vert T_{\Omega }\Vert
_{S^{2}}}{2}\text{,}
\end{equation*}%
where $T_{\Omega }$ depends on the window $g$ (or on the kernel of the
Toeplitz operator) and $N_{\Omega }$ is the largest integer less than $%
\left\vert \Omega \right\vert $ (this restriction is fundamental for the
methods used in \cite{APR}). 
\end{remark}

\begin{acknowledgement}
The author is thankful to Michael Speckbacher for bringing this problem to
his attention and for a wealth of discussions and collaborations about this
and related problems. Acknowledgements are also due to Federico Riccardi,
for stimulating discussions about the energy distribution of individual
eigenvalues, to Franz Luef, for helping to keep the author's interest in quantum
information theory alive, and to Fabio Nicola for encouragement and for
sharing the preprint \cite{NicolaRiccardiTilli2026}. GPT 5.6 Sol has been
used to help eliminate typos, mistakes and inconsistencies from the text.
This research was funded in part by the Austrian Science Fund (FWF) through
the project 10.55776/PAT8205923 (L.D.A.). For open access purposes, the
author has applied a CC BY public copyright license to any author-accepted
manuscript version arising from this submission.
\end{acknowledgement}

\end{document}